\documentclass[11pt]{article}
\usepackage[usenames]{color}
\usepackage{times,amsmath,amssymb,graphics,epsfig,ifthen}
\usepackage[colorlinks=true,
linkcolor=webgreen,
filecolor=webbrown,
citecolor=webgreen]{hyperref}

\definecolor{webgreen}{rgb}{0,.5,0}
\definecolor{webbrown}{rgb}{.6,0,0}

\newcommand\toPrint{TRANS}

\newtheorem{thm}{Theorem}[section]
\newtheorem{cor}[thm]{Corollary}
\newtheorem{lem}[thm]{Lemma}

\def\slfrac#1#2{\hbox{\kern.1em %
 \raise.5ex\hbox{\the\scriptfont0 #1}\kern-.11em %
 /\kern-.15em\lower.25ex\hbox{\the\scriptfont0 #2}}}

\DeclareMathOperator{\Trace}{Trace}

\newcommand{\al}{\alpha}
\newcommand{\QQ}{{\mathbb Q}}

\newcommand{\eeq}{\end{equation}}
\newcommand{\beql}[1]{\begin{equation}\label{#1}}
\newcommand{\eqn}[1]{(\ref{#1})}

\newcommand{\ZZ}{{\mathbb {Z}}}
 
\newcommand{\Complex}{{\rm {\bf C}\mkern-9mu\rule{0.05em}{1.4ex}\mkern10mu}}
\newcommand{\C}{{\rm {\bf C}\mkern-9mu\rule{0.05em}{1.4ex}\mkern10mu}}
\newcommand{\be}{\begin{equation}}
\newcommand{\ee}{\end{equation}}
\newcommand{\benonum}{\begin{displaymath}}
\newcommand{\eenonum}{\end{displaymath}}

\newcommand{\xbf}{\mbox{${\bf x }$} }
\newcommand{\ybf}{\mbox{${\bf y }$} }
\newcommand{\zbf}{\mbox{${\bf z }$} }

\newcommand{\Hbf}{\mbox{${\bf H }$} }
\newcommand{\Ibf}{\mbox{${\bf I }$} }

\newcommand{\Xbf}{\mbox{${\bf X }$} }
\newcommand{\Ybf}{\mbox{${\bf Y }$} }
\newcommand{\Zbf}{\mbox{${\bf Z }$} }

\makeatletter
\def\@sect#1#2#3#4#5#6[#7]#8{\ifnum #2>\c@secnumdepth
     \def\@svsec{}\else
     \refstepcounter{#1}\edef\@svsec{\csname the#1\endcsname.\hskip .75em }\fi
     \@tempskipa #5\relax
      \ifdim \@tempskipa>\z@
        \begingroup #6\relax
          \@hangfrom{\hskip #3\relax\@svsec}{\interlinepenalty \@M #8\par}%
        \endgroup
       \csname #1mark\endcsname{#7}\addcontentsline
         {toc}{#1}{\ifnum #2>\c@secnumdepth \else
                      \protect\numberline{\csname the#1\endcsname}\fi
                    #7}\else
        \def\@svsechd{#6\hskip #3\@svsec #8\csname #1mark\endcsname
                      {#7}\addcontentsline
                           {toc}{#1}{\ifnum #2>\c@secnumdepth \else
                             \protect\numberline{\csname the#1\endcsname}\fi
                       #7}}\fi
     \@xsect{#5}}
\def\@begintheorem#1#2{\it \trivlist \item[\hskip \labelsep{\bf #1\ #2.}]}
\makeatother

\newtheorem{example}{Example}[section]
\newtheorem{definition}{Definition}[section]

\newcommand{\Expt}{\mbox{${\mathbb E}$} }

\newcommand{\A}{\mathcal{A}}

\newcommand{\F}{GF}
\newcommand{\FF}{{\mathbb F}}

\newenvironment{proof}{
\noindent
{\bf Proof:}
}{
\hfill$\blacksquare$
}

\begin{document}
\begin{center}

\ifthenelse{\equal{\toPrint}{TRANS}}
{
{\Large{\bf Nonintersecting Subspaces Based on Finite Alphabets }} \\
}
{
{\Large{\bf Nonintersecting Subspaces Based on Finite Alphabets (Short Version) }} \\
}

\vspace{1\baselineskip}

\ifthenelse{\equal{\toPrint}{TRANS}}
{
{\em Fr\'ed\'erique E. Oggier}\footnote{This
work was carried out during F.~E.~Oggier's
visit to AT\&T Shannon Labs during the summer
of 2003. She thanks the
Fonds National Suisse,
Bourses et Programmes d'\'Echange for support.} \\
D\'epartement de Math\'ematiques  \\
Ecole Polytechnique F\'ed\'erale de Lausanne  \\
1015 Lausanne - Switzerland \\
(Email: frederique.oggier@epfl.ch) \\
\vspace{1\baselineskip}

{\em N. J. A. Sloane} \\
Information Sciences Research Center \\
AT\&T Shannon Labs \\
Florham Park, NJ 07932--0971, USA \\
(Email: njas@research.att.com) \\
\vspace{1\baselineskip}

{\em A. R. Calderbank} \\
Program in Applied and Computational Mathematics\\
Princeton University \\
Princeton, NJ 08540, USA \\
(Email: calderbank@math.princeton.edu) \\
\vspace{1\baselineskip}

{\em Suhas N. Diggavi} \\
School of Computer and Communication Sciences \\
Ecole Polytechnique F\'ed\'erale de Lausanne  \\
1015 Lausanne - Switzerland \\
(Email: suhas.diggavi@epfl.ch) \\
\vspace{2\baselineskip}
}
{
{\em Fr\'ed\'erique E. Oggier},
Ecole Polytechnique F\'ed\'erale de Lausanne  \\
\vspace{.5\baselineskip}

{\em N. J. A. Sloane},
AT\&T Shannon Labs \\
\vspace{.5\baselineskip}

{\em A. R. Calderbank},
Princeton University \\
\vspace{.5\baselineskip}

{\em Suhas N. Diggavi},
Ecole Polytechnique F\'ed\'erale de Lausanne  \\
\vspace{1\baselineskip}
}

{\bf Abstract} \\
\vspace{.5\baselineskip}
\end{center}


\setlength{\baselineskip}{1.5\baselineskip}
Two subspaces of a vector space are here called ``nonintersecting''
if they meet only in the zero vector.
The following problem arises in the design of noncoherent 
multiple-antenna communications systems.
How many pairwise nonintersecting $M_t$-dimensional subspaces 
of an $m$-dimensional vector space $V$ over a field $\FF$
can be found, if the generator matrices for the subspaces 
may contain only symbols from a given finite alphabet $\A \subseteq \FF$?
The most important case is when $\FF$ is the field of 
complex numbers $\Complex$; then $M_t$ is the number of antennas.
If $\A = \FF = \F(q)$ it is shown that 
the number of nonintersecting subspaces is
at most $(q^m-1)/(q^{M_t}-1)$, and that this bound can be attained
if and only if $m$ is divisible by $M_t$.
Furthermore these subspaces remain nonintersecting
when ``lifted'' to the complex field.
Thus the finite field case is essentially completely solved.
In the case when $\FF = \Complex$
only the case $M_t=2$ is considered. It is shown that if
$\A$ is a PSK-configuration, consisting of 
the $2^r$ complex roots of unity, the number of
nonintersecting planes is 
at least $2^{r(m-2)}$ and at most $2^{r(m-1)-1}$
(the lower bound may in fact be the best that can be achieved).

\section{Introduction}\label{sec:intro}

In \cite{Foschini96}, \cite{Telatar99} it was shown that the capacity of the
multiple-antenna channel grows linearly as
a function of the minimum of the numbers of
transmitting and receiving antennas.
The proof assumed that the receiver has complete information about
the channel.
In \cite{Tarokh98} the emphasis was placed on reducing error probability
by introducing correlation between signals transmitted from
different antennas.
These points of view can be combined by observing
that there is a trade-off between
rate and reliability \cite{Tarokh98}, \cite{ZhengTse02b}.

\ifthenelse{\equal{\toPrint}{TRANS}}
{
Most of the early work on multiple-antenna communications 
assumed that the receiver was able to track
the channel perfectly---i.e. used coherent detection.
If coherent detection is difficult or too expensive,
one can use noncoherent detection, as studied in 
\cite{HochwaldMarzetta99}.  The main result from this work
is that the capacity is still (almost) linear in the minimum number
of transmitting or receiving antennas \cite{HochwaldMarzetta99},
\cite{ZhengTse02a}. Hence, both in the coherent and
noncoherent cases, it was established that the use
of multiple-antennas leads to a gain in information transmission rate.
}{}

In \cite{HochwaldMarzetta00}, the
error probability of multiple-antenna noncoherent
communication channels was investigated.  It was
shown there (and in \cite{ZhengTse02a})
that if the channel is not known to the receiver,
the coding problem is equivalent to one of packing
subspaces (which represent codewords) according to a certain
notion of distance. 
The diversity order (the slope of
the error probability with respect to SNR) was shown to depend on the
dimension of the intersection of the subspaces.

In particular, to obtain maximal
diversity, one wishes to construct a family of
subspaces which intersect only at the origin. 
By a slight abuse of notation we will say that two
vector spaces are ``nonintersecting'' 
if their only common point is the zero vector.
A similar problem has been studied in
the context of designing differential codes for the multiple-antenna
channel \cite{HochSweld00}, \cite{Hughes00},
\cite{TarokhJaf00}. An extensive
characterization and classification of group differential space-time
codes was given in \cite{Shokrollahi01}. The focus of much of
this work is on constructing codes which have the nonintersecting
subspace property without imposing any constraints
on the  number of different symbols used to define
the codewords---that is, the codewords are allowed to use a
signal constellation that is larger than the minimum possible.

The main question addressed in the present paper 
is the construction of nonintersecting subspaces, subject
to the constraint that the codewords are defined using
symbols from a fixed, small constellation.
We focus on two cases: one in which
the symbols are taken from a finite field
and the other where they are taken from a PSK arrangement,
i.e. are complex roots of unity.
Our aim is to find
constructions that give the largest number of
nonintersecting subspaces (i.e. have
the highest rate) subject to these constraints.

It is worth remarking that a recent paper by Lusina et al.
\cite{LGB03} discusses an analogous problem 
for the case of coherent decoders.
Another related paper is Lu and Kumar
\cite{LK03} explores  code constructions with fixed alphabet
constraints for achieving different points on the rate-diversity trade-off.
Again, only coherent decoders are considered.
A very recent paper by Kammoun and Belfiore \cite{KaBe03}
directly addresses the problem of constructing codes
for non-coherent systems with a large value of
$\Lambda (\Xbf, {\Xbf}')$
(see \eqn{eq:PrinAng})
between subspaces.
However, their approach is quite different from ours.

The present paper is organized as
follows. In Section \ref{sec:chan}, we establish notation and formalize
the question being studied. In Section \ref{sec:FinField}, we study
the case when the symbols are taken from a finite field,
and in Section \ref{sec:PSK} when they are complex roots of unity
(i.e. PSK constellations).
Section \ref{sec:disc} compares the
different constructions and mentions some directions
for further research.

\section{Preliminaries}\label{sec:chan}
Let the number of transmitting antennas be $M_t$ and the number of receiving
antennas be $M_r$.
If $\ybf(k)\in\Complex^{M_r}$ is the received (column) vector at time $k$,
we can write
\begin{equation}\label{eq:model}
\ybf(k) = \sqrt{E_s} \, \Hbf(k) \, \xbf(k) + \zbf(k) ~,
\end{equation}
where the matrix $\Hbf(k)\in\Complex^{M_r\times M_t}$ represents the 
channel, the column vector $\xbf(k)\in\Complex^{M_t}$ is the channel input,
$E_s$ is the signal power per transmitting antenna, and
$\zbf(k)\in\Complex^{M_r}$ is zero mean i.i.d. Gaussian noise with
$\Expt[\zbf(k)\zbf(k)^H]=N_0\Ibf$.  We assume a
Rayleigh flat fading model, i.e. that the elements of
$\Hbf(k)$ are i.i.d. with a zero mean complex Gaussian distribution of unit variance.
The channel is assumed to be block time-invariant, that is,
$\Hbf(k)$ is independent of $k$ over a transmission block of $m$ symbols,
say $\Hbf(k) = \Hbf$
(although $\Hbf(k)$ may vary from block to block).
Looking at a single block of length $m$,
during which the channel is assumed to be time-invariant, we can
write
\begin{equation}\label{eq:model1}
\Ybf = [\ybf(1),\ldots,\ybf(m)] = \sqrt{E_s} \, \Hbf \, [\xbf(1),,\ldots,\xbf(m)] + 
[\zbf(1),\ldots,\zbf(m)] = \sqrt{E_s}\Hbf\Xbf + \Zbf ~.
\end{equation}
The focus of this paper is on constructing the space-time codewords $\Xbf$,
subject to
the constraint that the elements of $\Xbf$ are selected from a
particular alphabet $\A$.

\subsection{Criteria for code design}\label{subsec:PEP}

\ifthenelse{\equal{\toPrint}{TRANS}}
{
In this paper we assume that the receiver will not attempt to estimate the
channel matrix $\Hbf$, i.e. that we have a noncoherent receiver.
Therefore, the maximum likelihood 
detection rule without using the channel state information
(\cite{HochwaldMarzetta00}, \cite{Hughes00}) is that we
should decode $\Ybf$ as that codeword $\hat{\Xbf}$ which maximizes
\begin{equation}
\label{eq:detRule}
\frac{\exp(-\Trace[\Ybf \Psi ^{-1}\Ybf^H])}{|\pi\Psi|^{M_r}} ~,
\end{equation} 
where $\Psi =\Ibf+E_s\Xbf^H\Xbf$,
$H$ denotes the transposed complex conjugate or adjoint matrix,
and $| \cdot |$ denotes a determinant.
In the absence of channel state information at the receiver, 
Hochwald and Marzetta \cite{HochwaldMarzetta00}
argue that for high SNR, the
one should use unitary codewords $\Xbf$, satisfying
$\Xbf\Xbf^H=m\Ibf$.
Using this in (\ref{eq:detRule}) and the matrix inversion lemma 
(\cite[p. 19]{HJ85}),
it follows that $\hat{\Xbf}$ should be chosen to maximize
\begin{equation}
\label{eq:MLnoncoh}
\Trace[\Ybf\Xbf^H\Xbf\Ybf^H] ~.
\end{equation}
This implies that the decoder should
project the received signal onto the subspace
defined by each of the codewords
and declare the codeword with the maximal projection
to be the winner.
Using a Chernoff bound argument,
we find that the probability that a transmitted codeword $\Xbf$
is decoded as the codeword $\hat{\Xbf}$ is bounded above by
(\cite{HochwaldMarzetta00})
\begin{equation}
\label{eq:PEP}
\frac{1}
{|\Ibf_{M_t}+\frac{\rho^2 m^2}{4(1+\rho m)}[\Ibf_{M_t}-\frac{1}{m^2}
\hat{\Xbf}\Xbf^H\Xbf\hat{\Xbf}^H]|^{M_r}} ~,
\end{equation}
where $\rho=\frac{E_s}{N_0}$ is the signal-to-noise ratio (SNR).
If the SNR is large, this pairwise error probability
behaves like $(\frac{\Lambda\rho }{4})^{-M_r\nu}$, where $\nu$
is the rank of
$[\Ibf_{M_t}-\frac{1}{m^2}\hat{\Xbf}\Xbf^H\Xbf\hat{\Xbf}]$,
$$
\Lambda = \Lambda(\Xbf, \hat{\Xbf}) =
|m\Ibf_{M_t}-\frac{1}{m}\hat{\Xbf}\Xbf^H\Xbf\hat{\Xbf}^H|_+
^{\frac{1}{\nu}} ~,
$$
and $|\cdot|_+$ denotes the product of the nonzero
eigenvalues. Note that
$$
\left |\left[\begin{array}{c}\Xbf \\ \hat{\Xbf}\end{array}\right ] \left [
\begin{array}{cc}\Xbf^H & \hat{\Xbf}^H \end{array}\right ]\right | = 
|m^2\Ibf_{M_t}-\hat{\Xbf}\Xbf^H\Xbf\hat{\Xbf}^H| ~,
$$
which shows that $\nu=M_t$ is
equivalent to the condition that the
rows of $\Xbf,\hat{\Xbf}$ are linearly independent
(\cite{Hughes00}).
For this to happen we must have $m\geq 2M_t$.
}{}

Another interpretation can be given
in terms of the principal angles between subspaces
corresponding to pairs of codewords.
The principal angles between subspaces $\Xbf$ and $\Xbf '$ are given by
$\cos \theta_i =\frac{1}{m}\sigma_i({\Xbf}' \Xbf^H)$ where
$\sigma_i(\cdot)$
is the $i$-th singular value of the matrix (\cite{grass1}, \cite{GVL89}).
Using this we obtain
\begin{equation}\label{eq:PrinAng}
\Lambda (\Xbf, {\Xbf}') = m\prod_{i=1}^{\nu}[1-\cos^2 \theta_i]
=m\prod_{i=1}^{\nu}\sin^2 \theta_i ~.
\end{equation}
This provides a better measure of how good a code is:
not only should the subspaces be nonintersecting,
the value of $\Lambda (\Xbf, {\Xbf}')$
should be large for every pair $\Xbf$, $\Xbf '$ of distinct subspaces.
The error probability will be dominated by the pair
of codewords
with the least rank $\nu$ and the least ``distance''
$\Lambda (\Xbf, {\Xbf}')$.
For well separated subspaces this ``distance'' can also be approximated
by 
\begin{equation}\label{eq:PrinAng2}
\sum_{i=1}^{\nu}\sin^2 \theta_i ~,
\end{equation}
which is the
the notion of distance between subspaces used in \cite{grass1}
and \cite{grass3}.

Another way to compare these codes is by using the 
notion of diversity order (cf. \cite{Tarokh98}).

\begin{definition}
\label{defn:div}
If the average error probability
$\bar{P}_e(\mbox{$\rho$})$ as a function of the SNR $\rho$ satisfies
\begin{equation}
\label{eq:DivDefn}
\lim_{\mbox{$\rho$}\rightarrow \infty}
\frac{\log(\bar{P}_e(\mbox{$\rho$}))}{\log(\mbox{$\rho$})} = -d ~,
\end{equation}
the coding scheme is said to have {\em diversity order} $d$.
\end{definition}
It follows from (\ref{eq:PEP})
that the diversity order of the coding scheme
is equal to $M_r\nu$. The maximal diversity order
that can be achieved is therefore $M_rM_t$. We call codes that achieve
this bound {\em fully diverse} codes.

In brief, to get a diversity order of $M_rM_t$,
we need to construct nonintersecting subspaces which are far apart
in the metric defined by (\ref{eq:PrinAng}).
In this paper we will focus on obtaining maximal diversity
order by constructing families of
subspaces which are nonintersecting. In order to further
improve performance we need to maximize $\Lambda (\Xbf, {\Xbf}')$ over all
pairs $\Xbf$, $\Xbf '$ of distinct subspaces.
The rate of a code $C$ is 
$R = \frac{1}{m}\log(|\mathcal{C}|)$.
In trying to construct the maximal number of non-intersecting subspaces, we attempt
to get the highest rate codes that achieve maximal diversity order.

\subsection{Statement of the problem}\label{subsec:prelim}

\begin{definition}\label{def:disjoint}
Let $\FF$ be a field.
A {\em codeword} or {\em subspace} will mean an $M_t$-dimensional subspace
of $\FF^m$.  Two subspaces $\Pi _1$ and $\Pi _2$
are said to be {\em nonintersecting} over $\FF$ if 
their intersection is trivial, i.e. if
$ \Pi_1 \cap \Pi_2 = \{0\}$.
\end{definition}

\ifthenelse{\equal{\toPrint}{TRANS}}
{
Suppose $\Pi _1$ is generated by (row) vectors
$u_1, \ldots, u_{M_t} \in \FF^m$, and
$\Pi _2$ is generated by vectors
$v_1, \ldots, v_{M_t} \in \FF^m$.
Let 
$P:=\left[ \begin{array}{c} \Pi_1 \\ \Pi_2 \end{array} \right]$
denote the $2M_t \times m$  matrix with rows
$u_1, \ldots, u_{M_t}, v_1, \ldots, v_{M_t}$.
Then the following lemma is readily established.

\begin{lem}\label{lemma1}
The following properties are equivalent:
(i) $\Pi_1$ and $\Pi_2$ are nonintersecting,
(ii) $P$ has rank $2M_t$ over $\FF$, and
(iii) if $m = 2M_t$ the determinant of $P$ is nonzero.
\end{lem}
}{}

Suppose now that instead of 
allowing the entries in the matrices $\Pi _1$ and $\Pi _2$ to
be arbitrary elements of $\FF$, we restrict them
to belong to a finite subset
$\A \subseteq \FF$, called the {\em alphabet}.
In other words, the vectors $u_1, \ldots, u_{M_t}, v_1, \ldots, v_{M_t}$
must belong to  $\A^m$.
The question that we address is the following: given $M_t$,
$m$ and a finite alphabet $\A \subseteq \FF $, how many subspaces
can we find which are generated by vectors from
$\A^m$ and which are pairwise nonintersecting over $\FF$?
Furthermore, if the size of $\A$ is specified in advance,
which choice of $\A$ permits the biggest codes?

\ifthenelse{\equal{\toPrint}{TRANS}}
{
We first dispose of the trivial case when $M_t = 1$.
Two nonzero vectors $u, v$ are said to be
{\em projectively distinct} over a field $\FF$ if there is no $a \in F$
such that $u = av$.
Then if $M_t=1$, the maximum number of nonintersecting
subspaces is simply the maximum number of 
projectively distinct vectors in $\A ^m$.
}{}

In the following sections we will investigate
the first question for two kinds of alphabets:
(a) $\A$ is a finite field $\FF$ (Section \ref{sec:FinField}), and
(b) $M_t=2$ and  $\A \subseteq \C^m$ is a set of complex roots of unity
(Section \ref{sec:PSK}).

\ifthenelse{\equal{\toPrint}{TRANS}}
{
Of course, for the application to multiple-antenna code design,
the subspaces need to be disjoint over $\Complex$. 
In Theorem \ref{th:lift} of  Section \ref{sec:FinField}
we translate the results obtained over 
$\FF $ to this case by ``lifting'' the subspaces to the complex
field.
Furthermore, for this application, the case 
$m = 2 M_t$ is the most important.
}{}

\section{Finite Fields}\label{sec:FinField}

In this section we assume that 
the alphabet $\A$ and the field $\FF$ are both
equal to the finite field $\F(q)$, where 
$q$ is a power of a prime $p$.
At the end of the section we show 
how to ``lift'' these planes to the complex field
(see Theorem \ref{th:lift}).
In this case there is an obvious upper bound which can be
achieved in infinitely many cases.
Let $V$ denote the vector space $\F(q)^m$.

\begin{thm}\label{prop:bnd1}
The number of pairwise nonintersecting $M_t$-dimensional
subspaces of $V$ is at most
\beql{Eq1}
\frac{q^m -1}{q^{M_t}-1}  ~.
\eeq
\end{thm}

\begin{proof}
There are $q^m-1$ nonzero vectors in $V$ and each subspace
contains $q^{M_t}-1$ of them.
No nonzero vector can appear in more than one subspace.
\end{proof}

It is convenient here to use the language of projective geometry,
c.f. \cite[Appendix B]{MS77}.  
Recall that the points of the projective space
$P(s,q)$ are equivalence classes of nonzero vectors from $GF(q)^{s+1}$,
where two vectors are regarded as equivalent if one is
a nonzero scalar multiple of the other.

A {\em spread}
\cite{Hirsch} in $PG(s,q)$ is a partition of the
points into copies of $PG(r,q)$.

\begin{thm}
Such a spread exists
if and only if $r+1$ divides $s+1$.
\end{thm} 

\begin{proof}
This is a classical result, due to Andr\'e
(\cite{Andre}; \cite[Theorem~4.1.1]{Hirsch}).
\end{proof}

\begin{cor}\label{cor:FF}
The bound \eqref{Eq1} can be attained
whenever $M_t$ divides $m$, and only in those cases.
\end{cor} 

\begin{proof}
This is immediate from the theorem, since a set of
points in a projective space represents
a set of projectively distinct
lines in the corresponding vector space.
\end{proof}

Note that the condition is independent of $q$.
If a set of nonintersecting subspaces meeting
\eqref{Eq1} exists over one finite field then 
it exists over every finite field.


Furthermore, it is straightforward to
construct the nonintersecting subspaces 
meeting the bound in \eqref{Eq1},
as we now show.
The nonzero elements of a finite field
$\FF$ form a multiplicative group which will be
denoted by $\FF^{\ast}$. This is a cyclic
group \cite[Chap. 2]{LN83}.

Suppose $M_t$ divides $m$, and consider
the fields $F_0 = \F(q)$, 
$F_1 = \F(q^{M_t})$, $F_2 = \F(q^m)$.
Then $F_0 \subseteq F_1 \subseteq F_2$. 
By regarding $\F(q^m)$ as a vector space of dimension $m$ over 
$\F(q)$ we can identify $F_2$ with $V$.
Similarly we can regard $F_1$ as a
$M_t$-dimensional subspace of $V$.
The desired spread is now obtained by partitioning $F_2^{\ast}$
into (multiplicative) cosets of $F_1^{\ast}$.

\begin{example} \rm
We consider the case $M_t=2$, $m=4$ and $\A = \F(2) = \{0,1\}$.
Then $F_0 = \F(2), F_1 = \F(4), F_2 = \F(16)$.
Each plane in $\F(2) ^4$ contains three nonzero vectors, and $\F(2)^4$
itself contains 15 nonzero vectors.
We wish to find a spread of $PG(1,2)$'s inside $PG(3,2)$,
that is, a partitioning of the 15 vectors into five disjoint sets of
three, where each set of three adds to
the zero vector.

Let $\F(16) = \F(2)[\al]$, where $\al ^4 + \al + 1 = 0$.
A table of the elements of this field and their
binary representations can be found for
example in \cite[Fig. 3.3]{MS77}.
Then $\F(4)$ is the subfield $\{1, \al ^5, \al^{10}\}$,
so $F_1 ^{\ast}  = \{\al ^5, \al^{10}\}$, and we obtain
the desired partition 
$$
F_2^{\ast} = \bigcup_{j=0}^{4} \al^j F_1^{\ast} ~.
$$
Only two of the three vectors are needed to define each
plane, so we have the following generators for
the five planes:
$$
(1,\al), ~(\al,\al^6), ~(\al^2, \al^7), ~(\al^3, \al^8), ~(\al^4, \al^9) ~.
$$
Using the table in \cite{MS77}, we convert these to explicit 
generator matrices for the five nonintersecting planes:
$$
\left[ \begin{array}{l} 1000 \\
                        0110 \\
              \end{array} \right] ,
\left[ \begin{array}{l} 0100 \\
                        0011 \\
              \end{array} \right] ,
\left[ \begin{array}{l} 0010 \\
                        1101 \\
              \end{array} \right] ,
\left[ \begin{array}{l} 0001 \\
                        1010 \\
              \end{array} \right] ,
\left[ \begin{array}{l} 1100 \\
                        0101 \\
              \end{array} \right] .
$$
\end{example}

The problem is therefore essentially solved as long as
$M_t$ divides $m$.  If not, we can use
partial spreads--see the surveys in
\cite{Eisfeld} and \cite{Soicher1}.

We end this section by observing
that a set of nonintersecting subspaces over
a finite field $\A = \F(q)$, $q=p^k$, $p$ prime, can always be ``lifted'' to
a set of nonintersecting subspaces over a complex
alphabet $\bar{\A}$ of the same size.

This can be done as follows.
Suppose $\F(q)= \F(p)[\al]$, where $\al$ is a root
of a primitive irreducible polynomial $f(X) \in \F(p)[X]$.
Let $n = p^k-1$ and let 
$\mu _n = e^{2 \pi i / n}$.
Adjoining $\mu_n$ to the rational numbers $\QQ$, we obtain
the cyclotomic field $\QQ(\mu_n)$, with ring of
integers $\ZZ[\mu_n]$.
It is a classical result from number theory that
the ideal $(p)$ in $\ZZ[\mu_n]$ factors into
$g = \varphi (n)/k$ distinct maximal prime ideals
$\mathfrak{p}_1, \mathfrak{p}_2, \ldots, \mathfrak{p}_g$,
where $\varphi(\cdot)$ is the Euler totient function.
Furthermore, for each $\mathfrak{p}_j$,
the residue class ring $\ZZ[\mu_n]/\mathfrak{p}_j  \cong \F(q)$
(see for example
\cite[Theorem 10.45]{Cohn},
\cite[Chap. 10, \S 3B]{Rib},
\cite[Theorem 2.13]{Wash},
\cite[Theorem 7-2-4]{Weiss}).
If we choose $\mathfrak{p}_j$ to be the ideal generated by $p$ and
$f(\mu_n)$, then $\ZZ[\mu_n]/\mathfrak{p}_j$
is exactly the version of $GF(q)$ that we started with.
Note that since
$\mathfrak{p_j}$ contains $(p)$, it acts as reduction mod $p$ on $\ZZ$.
We therefore have a ring
homomorphism from $\ZZ[\mu_n]$ to $\F(q)$ given by
\beql{EqLift}
\phi : \ZZ[\mu_n] \stackrel{{\rm mod~} p}{\to} \ZZ[\mu_n]/\mathfrak{p_j} \stackrel{\cong}{\to} \F(q) ~.
\eeq
In this way we can lift vectors over $GF(q)$ to
vectors over the alphabet $\bar{A}$ consisting of $0$ and the 
$q-1$ powers of $\mu_n$.

%

Example: Let $GF(8) = GF(2)[\alpha]$
where $\alpha$ is a root of $X^3+X+1$.  Then $q=8$, $n=7$, 
$\mu_7 = e^{2 \pi i/7}$.
To lift $GF(8)$ to $\Complex$ we write
$GF(8) = \{ 0, 1, \alpha, \alpha^2, \ldots, \alpha^6 \}$,
and lift $0$ to $0$ and $\alpha^j$ to $\mu_7^j$ for $j=0, \ldots, 6$.

Let $\Pi$ be an $M_t$-dimensional subspace of $\F(q)^m$. By lifting each
element of a generator matrix we obtain an $M_t$-dimensional subspace
$\bar{\Pi} \subseteq \Complex^m$, defined over
an alphabet $\bar{\A}$ of size $q$.

\begin{thm}\label{th:lift}
If two subspaces $\Pi_1, \Pi_2$ of $\F(q)^m$ are
nonintersecting, so are their lifts $\bar{\Pi}_1, \bar{\Pi}_2$.
\end{thm}

\begin{proof}
Let 
$P:=\left[ \begin{array}{c} \Pi_1 \\ \Pi_2 \end{array} \right]$
and
$\bar{P}:=\left[ \begin{array}{c} \bar{\Pi_1} \\ \bar{\Pi_2} \end{array} \right]$.
By Lemma \ref{lemma1}, $P$ has a $2M_t \times 2M_t$ invertible
submatrix. Since $\phi$ is a ring homomorphism, the
lift of this submatrix is also invertible.
\end{proof}

It follows that the subspaces constructed in
Corollary \ref{cor:FF} are also nonintersecting
when lifted to the complex field.

This construction gives full diversity order non-coherent
space-time codes when
the elements of the codewords are restricted to belong to a finite field.
Their rate is
\begin{displaymath}
R = \frac{1}{m}\log(q^m-1)-\frac{1}{m}\log(q^{M_t}-1) < \log(q) ~,
\end{displaymath}
which according to Theorem \ref{prop:bnd1} is the maximal achievable
rate for diversity order $M_tM_r$. Moreover, the above
relationship implies that for fully diverse codes constructed
from a finite field, we cannot achieve a rate higher than $\log(|\A|)$.


\section{PSK constellations}\label{sec:PSK}

Throughout this section we assume that the alphabet $\A$ consists of
the set of complex $2^r$-th roots of unity, 
that is, $\A = \{e^{2 \pi i j / 2^r}, 0 \le j < 2^r \}$, for some
$r \ge 1$.
Let $\mu = e^{2 \pi i / 2^r}$ be a primitive $2^r$-th root of unity;
$\A$ is a cyclic multiplicative group with generator $\mu$.
In this section we assume that $M_t=2$, that is,
the code consists of a set of pairwise
nonintersecting planes.

\ifthenelse{\equal{\toPrint}{TRANS}}
{
\begin{example}\label{ex:unity} \rm
Some examples of roots of unity:
\begin{enumerate}
\item If $r=1$, $\mu=-1$ and the alphabet is $\A=\{1,-1\}$.
\item If $r=2$, $\mu= i$ and the alphabet is $\A=\{1,i,-1,-i \}$.
\item If $r=3$, $\mu = (1+i)/\sqrt{2}$ and the alphabet 
is $\A = \{e^{\pi i j / 4}, 0 \le j \le 7 \}$.  This is the 8-PSK constellation. 
\end{enumerate}
\end{example}
}{}

There is a trivial upper bound.
\begin{thm}\label{Th41}
Let $\A$ be the set of $2^r$ roots of unity, $r \geq 1$. 
Then the number of pairwise nonintersecting planes is at most
$\frac{1}{2} |\A|^{m-1} = 2^{(m-1)r-1}$.
\end{thm}

\ifthenelse{\equal{\toPrint}{TRANS}}
{
\begin{proof}
If $v_1, v_2 \in \A^m$ are the generators for a plane,
that plane also contains all multiples $\mu^j v_1$ and $\mu^j v_2$,
a total of $2 |\A|$ vectors. Since
these sets of vectors must all be disjoint, the number of planes is at most 
$|\A|^m / (2 |\A|)$.
\end{proof}

The same argument shows that
there are at most $\frac{1}{M_t} |\A|^{m-1}$
nonintersecting $M_t$-dimensional subspaces of complex
$m$-dimensional space for any finite alphabet $A$.
The implication of this in
terms of rate is that
\begin{displaymath}
R \leq \frac{m-1}{m}\log(|\A|)- \frac{1}{m}\log(M_t) < \log(|\A|) ~.
\end{displaymath}
Hence, for fully diverse codes constructed from PSK constellations, we
cannot achieve a rate exceeding $\log(|\A|)$.

\begin{example}\label{ex:bound} \rm
Let $\A$ be the set $\{1,i,-1,-i\}$ and take $m=4$. The total number 
of vectors in $\A^4$ is $4^4$. Each vector has 4 multiples, so 
each plane accounts for at least 8 vectors. 
Therefore there are at most $\frac{4^4}{8}=32$ planes.
\end{example}
}{}

In the other direction we will prove:

\begin{thm}\label{Th42}
Assume $r \ge 1$ and that $m \ge 2$ is even.
There exist $N = |\A|^{m-2} = 2^{(m-2)r}$ pairwise nonintersecting planes
in $\Complex^m$ defined using the complex $2^r$-th roots of unity.
\end{thm}

Note that the upper and lower bounds coincide 
in the case $r=1$, that is, when $\A = \{1, -1\}$.

\ifthenelse{\equal{\toPrint}{TRANS}}
{
The proof is simplified by the use of valuations
(cf. \cite{Gouvea}). 
If $x \in \QQ$, $x = 2^a \, \frac{b}{c}$ with
$a, b, c \in \ZZ, c \neq 0$, $b$ and $c$ odd, then 
the $2$-adic valuation of $x$ is $\nu_2(x) = a$.
Similarly, suppose $x$ belongs to the cyclotomic field $\QQ(\mu)$.
Since $1-\mu$ is a prime in $\ZZ[\mu]$, we can write $x$ uniquely as
$(1-\mu)^a \, \frac{b}{c}$ with
$a \in \ZZ$, $b, c \in \ZZ[\mu], c \neq 0$, $b$ and $c$ relatively prime to $1-\mu$.
The $(1-\mu)$-adic valuation of $x$ is then
$\nu_{1-\mu}(x) = a$.
It is easy to check that
for $k \in \ZZ$, $k \neq 0$,
$\nu_{1-\mu}(1-\mu^k)=2^{\nu_2(k)}$. In particular,
if $k \in \ZZ$ is odd, $\nu_{1-\mu}(1-\mu^k)=1$.

We will also need a lemma:

\begin{lem}\label{lemma4}
Let $\Pi$ be a plane in $\Complex^m$ generated by
vectors $v_1$, $v_2$, and denote by 
\[
\tilde{\Pi}_1=\left[\begin{array}{ccc}
                    v_1 & x_{11} & x_{12} \\
                    v_2 & x_{21} & x_{22} \\
                    \end{array} \right] 
\]
and 
\[
\tilde{\Pi}_2= \left[\begin{array}{ccc}
                    v_1 & y_{11} & y_{12} \\
                    v_2 & y_{21} & y_{22} \\
                    \end{array} \right]
\]
two different embeddings of $\Pi$ into $\Complex^{m+2}$. 
Then $\tilde{\Pi}_1 \cap \tilde{\Pi}_2 = \{0\}$ 
if and only if
\[
\bigg| \begin{array}{cc}
                 y_{11} - x_{11} & y_{12} - x_{12} \\
                  y_{21} - x_{21} & y_{22} - x_{22} \\
                 \end{array} \bigg| \neq 0.
\]
\end{lem}

\begin{proof}
By Lemma \ref{lemma1}, it is necessary and sufficient that
the matrix
$P:=\left[ \begin{array}{c} \tilde{\Pi}_1 \\ \tilde{\Pi}_2 \end{array} \right]$
have rank 4.
Subtracting the first and second rows of $P$
from the third and fourth rows, we get the matrix
\[
\left[ \begin{array}{ccc}
                v_1 & x_{11} & x_{12} \\
                v_2 & x_{21} & x_{22} \\
                0 & y_{11}-x_{11} & y_{12}-x_{12} \\
                0 & y_{21}-x_{21} & y_{22}-x_{22} \\
                \end{array} \right ] ~.
\]
and the result follows.
\end{proof}

We now give the proof of the theorem, for which we use induction
on even values of $m$.  For $m=2$ we take the single plane
\[
\left[ \begin{array}{rr}
                1 & 1 \\
                1 & -1 \\
                \end{array} \right ] ~.
\]
Suppose the result is true for $m$.  For each of the $|\A|^{m-2}$
pairwise nonintersecting
planes in $\Complex^m$ we will construct $|\A|^2$ 
planes in $\Complex^{m+2}$, such
that full set of planes so obtained is
pairwise nonintersecting; this will establish
the desired result.

If two planes are nonintersecting
in $\Complex^m$ then they are certainly nonintersecting
when embedded in any way in $\Complex^{m+2}$.
So we need only show that the $|\A|^2$ embeddings of any
single plane are pairwise nonintersecting.

Let $\Pi$ be a plane in $\Complex^m$ generated by
vectors $v_1$, $v_2$, and denote by $\tilde{\Pi}(a,b)$
the plane in $\Complex^{m+2}$ with generator matrix
\[
\left[\begin{array}{ccc}
                    v_1 & \mu ^a & \mu ^b \\
                    v_2 & \mu ^{a+b} & \mu ^{a+2b+1} \\
                    \end{array} \right]  ~,
\]
for $a, b = 0, 1, \ldots, 2^r - 1$.


We will use Lemma \ref{lemma4} to show that all the planes
$\{\tilde{\Pi}(a,b) \mid a \in \A, \, b \in \A \}$ are
pairwise nonintersecting. For this we must show that

\[
\bigg| \begin{array}{cc}
\mu^c - \mu^a & \mu^d - \mu^b \\ 
\mu^{c+d}-\mu^{a+b} & \mu^{c+2d+1}-\mu^{a+2b+1}   \\
\end{array} \bigg| ~=~ 0 
\]
if and only if $a = c$ and $b = d$.

The above determinant is equal to
\beql{Eq13}
\mu^{2c+2d+1}(1-\mu^{a-c})(1-\mu^{(a-c)+2(b-d)})
-\mu^{c+2d}(1-\mu^{b-d})(1-\mu^{(a-c)+(b-d)}) ~.
\eeq
If the determinant is zero, 
the $(1-\mu)$-adic valuations of the two
terms on the right must be equal, that is,
\begin{equation}\label{eq:valuation}
2^{\nu_2(a-c)}+2^{\nu_2(a-c+2(b-d))}= 2^{\nu_2(b-d)}+2^{\nu_2(a-c+b-d)} ~.
\end{equation}

We must show that this is true if and only if $a=c$ and $b=d$.
We consider four cases, depending on the parity of $a-c$ and $b-d$.
If  $a-c\equiv 1, b-d\equiv 1 \mbox{(mod~}2)$ then
\eqref{Eq13} reads
$1+1 = 1 + 2^{\nu_2(a-c+b-d)} \ge 3$
(since $a-c+b-d$ is even),
a contradiction.
Similarly,
if  $a-c\equiv 1, b-d\equiv 0 \mbox{(mod~}2)$ we get
$1+1 = 2^{\nu_2(b-d)}+ 1$,
and 
if  $a-c\equiv 0, b-d\equiv 1 \mbox{(mod~}2)$ we get
$2^{\nu_2(a-c)}+2^{\nu_2(a-c+2(b-d))} = 1 + 1$,
which are also contradictions. The fourth possibility is
$a-c \equiv b-d \equiv 0$ (mod 2).
Let $a-c = 2^s x$ and $b-d = 2^t y$, where $x$ and $y$ are odd, $s,t \geq 1$.
We have
\[
\nu_2(a-c+2(b-d))=
\left\{
\begin{array}{ll}
s & \mbox{if } s < t \\
s & \mbox{if } s= t \\
\geq t & \mbox{if } s>t 
\end{array}
\right.
\]
and
\[
\nu_2(a-c+2(b-d))=
\left\{
\begin{array}{ll}
s & \mbox{if } s < t \\
\geq s & \mbox{if } s= t \\
 t & \mbox{if } s>t 
\end{array}
\right.
\]
Substituting these valuations in equation (\ref{eq:valuation}) again
gives a contradiction. This concludes the proof of Theorem \ref{Th42}.
}{}

\section{Discussion}\label{sec:disc}

The following table compares the codes constructed in
Sections 3 and 4 in the case $M_t=2$,
i.e. codes which are pairwise nonintersecting
$2$-dimensional subspaces of $\Complex^m$,
for $m=4, 6$ and $8$, and
alphabets $\A$ of sizes 2, 4 and 8.
The top entry in each cell
gives the number of planes obtained from the finite field
construction (Corollary \ref{cor:FF}).
The bottom entry gives the lower and upper 
bounds obtained using complex $|\A|$-th roots of unity,
from Theorem \ref{Th42} and Theorem \ref{Th41}.
Asymptotically, the rates of the two constructions 
are very similar.
Both satisfy
$\log (\mbox{~number~of~codewords~}) / m \approx \log (|\A|)$,
for $m$ large, and so both asymptotically achieve the maximal rate
possible for fully diverse codes.

Note that the construction via finite fields
results in codes for which alphabet consists of
$0$ and the complex $(|A|-1)$-st roots of unity,
whereas the construction via PSK
constellations produces codes in which the symbols are
the complex $|A|$-th roots of unity (and $0$ is not used).

\vspace{1\baselineskip}

\begin{center}
$
\begin{array}{|c|c|c|c|}
\hline
       & m=4 & m=6 & m=8 \\
\hline
|\A|=2 & 5 & 21 & 85 \\
       & 4-4 & 16-16 & 64-64 \\
\hline
|\A|=4 & 17 & 273 & 4369 \\
       & 16-32 & 256-512 & 4096-8192 \\
\hline
|\A|=8 & 65 & 4161 & 266305 \\
       & 64-256 & 4096-16384 & 262144-1048576 \\
\hline
\end{array}
$
\end{center}
\begin{center}
Table I. Number of pairwise nonintersecting planes
in $\Complex^m$ for various \\
sizes of the alphabet $|\A|$  
(see text for details).
\end{center}

\vspace{1\baselineskip}

We end by mentioning some
topics for further research.

\begin{itemize}
\item
We also used clique-finding algorithms to
search for larger sets of planes
than those given in Theorem \ref{Th42},
again taking $\A$ to be the set of
$2^r$-th complex roots of unity.
These searches were unsuccessful,
and so we have not mentioned them elsewhere in the paper.
These negative results lead us to conjecture, albeit weakly,
that the lower bounds in Theorem \ref{Th42} cannot be improved.
It would be nice to have a better upper bound than
that in Theorem \ref{Th41} for the case $r > 1$.
It would also be a worthwhile project to do a more extensive
computer search for better codes, both for the above
alphabet and for other alphabets.

It is straightforward to formulate the search as 
a clique-finding problem.
The first step is to prepare a list of candidate
subspaces, making sure that the generator
matrices use only symbols from $\A$,
and that the subspaces have the specified dimension and are
distinct (a subspace may have many different generator matrices:
only one version is placed on the list of candidates).
Then a graph is constructed with the 
candidate subspaces as vertices, and with an edge
joining two vertices if and only if
the subspaces are nonintersecting.
Then a good code is a maximal clique
in this graph.


\item
Can the construction in Theorem \ref{Th42}
be generalized to the case when $M_t$ is
larger than $2$?
In particular, it would be interesting to do a computer search 
in the case $M_t=3$ and $m=6$.

\item 
This paper has focused only on the existence and construction of
finite alphabet codes which achieve maximal diversity order,
and we did not consider decoding complexity.
The decoding problem
involves projecting the received matrix $\Ybf$ onto the
candidate subspaces (see (\ref{eq:MLnoncoh})).
In general this may require a search over 
$2^{mR}$ codewords, where $R$ is the rate of the code. Since this
number grows exponentially with the code length,
a natural question to ask is
whether there are codes which are optimally decodable in polynomial time,
or have polynomial time sub-optimal decoders which
perform satisfactorily.

\item
In \cite{grass1} (see also \cite{grass3}) a large number
of optimal or putatively optimal
packings of subspaces in $\Complex^m$
were constructed using (\ref{eq:PrinAng2})
as a measure of ``distance'' between subspaces.
It would be worthwhile repeating these calculations using
(\ref{eq:PrinAng}) instead.

\end{itemize}

\end{document}